\documentclass[12pt,a4paper,reqno]{amsart}
\usepackage{amsmath}
\usepackage{amsthm}
\usepackage{amsfonts}
\usepackage{amssymb}
\usepackage{color}
\usepackage{mathrsfs}
\usepackage{graphicx}

\numberwithin{equation}{section}

\theoremstyle{plain}
\newtheorem{thm}{Theorem}

\theoremstyle{plain}
\newtheorem{lem}{Lemma}[section]

\theoremstyle{plain}
\newtheorem*{theorema}{Theorem~A}

\theoremstyle{remark}
\newtheorem*{remark}{Remark}

\def\bfa{\mathbf{a}}

\def\bfp{\mathbf{p}}

\def\bfv{\mathbf{v}}
\def\bfw{\mathbf{w}}

\def\dd{\mathrm{d}}
\def\ee{\mathrm{e}}
\def\ii{\mathrm{i}}

\def\eps{\varepsilon}

\def\Rr{\mathbb{R}}

\def\Zz{\mathbb{Z}}

\def\HHH{\mathcal{H}}
\def\III{\mathcal{I}}

\def\LLL{\mathcal{L}}
\def\MMM{\mathcal{M}}

\def\PPP{\mathcal{P}}

\def\SSS{\mathcal{S}}

\def\frakB{\mathfrak{B}}

\def\frakI{\mathfrak{I}}

\DeclareMathOperator{\meas}{meas}
\DeclareMathOperator{\area}{area}

\DeclareMathOperator{\rp}{RP}

\renewcommand{\le}{\leqslant}
\renewcommand{\ge}{\geqslant}

\title[Billiards in polyhedra]
{Billiards in polyhedra: a method\\
to convert 2-dimensional uniformity\\
to 3-dimensional uniformity}

\author[Beck]{J. Beck}

\address{Department of Mathematics, Hill Center for the Mathematical Sciences, Rutgers University, Piscataway NJ 08854, USA}

\email{jbeck@math.rutgers.edu}

\author[Chen]{W.W.L. Chen}

\address{School of Mathematical and Physical Sciences, Faculty of Science and Engineering, Macquarie University, Sydney NSW 2109, Australia}

\email{william.chen@mq.edu.au}

\author[Yang]{Y. Yang}

\address{School of Science, Beijing University of Posts and Telecommunications, Beijing 100876, China}

\email{yangyx@bupt.edu.cn}



\begin{document}

\begin{abstract}
The class of $2$-dimensional non-integrable flat dynamical systems has a rather extensive literature with many deep results,
but the methods developed for this type of problems, both the traditional approach via Teichm\"{u}ller geometry
and our recent shortline-ancestor method, appear to be exclusively plane-specific.
Thus we know very little of any real significance concerning $3$-dimensional systems.

Our purpose here is to describe some very limited extensions of uniformity in $2$ dimensions to uniformity in $3$ dimensions.
We consider a $3$-manifold which is the cartesian product of the regular octagonal surface with the unit torus.
This is a restricted system, in the sense that one of the directions is integrable.
However, this restriction also allows us to make use of a transference theorem for arithmetic progressions
established earlier by Beck, Donders and Yang.
\end{abstract}

\maketitle

\thispagestyle{empty}

%
%

\section{Introduction}\label{sec1}

For rational polygons where every angle is a rational multiple of~$\pi$, we have the following fundamental result
of Kerckhoff, Masur and Smillie~\cite{KMS86} in 1986.

\begin{theorema}
Let $P$ be a rational polygon.
For almost every initial direction and for every non-pathological starting point for this direction,
the half-infinite billiard orbit in $P$ is uniformly distributed.
\end{theorema}

Given any initial direction, a point $\bfp_0\in P$ is called a pathological starting point for this direction if the half-infinite billiard orbit
starting from $\bfp_0$ and with this direction hits a singularity of~$P$.
Otherwise the point $\bfp_0$ is called a non-pathological starting point for this direction.
It is easy to see that for any given direction, almost every point in $P$ is a non-pathological starting point.

The proof of Theorem~A consists of essentially three steps.
The first step is to establish the ergodicity of the corresponding interval exchange transformation.
The second step is to use the well known Birkhoff ergodic theorem.
The final step is to extend ergodicity to unique ergodicity.

Our aim is to convert Theorem~A to a result concerning equidistribution of $3$-dimensional billiard in some polyhedra.
However, we need to restrict our discussion to rational polygonal right prisms.
A rational polygonal right prism is a region in $3$-dimensional cartesian space of the form
\begin{equation}\label{eq1.1}
M=P\times I=\{(x,y,z)\in\Rr^3:(x,y)\in P\mbox{ and }z\in I\},
\end{equation}
where $P$ is a rational polygon and $I=[0,z_0]$ is an interval.

As the rational polygonal right prism $M=P\times I$ is integrable in the direction of the interval~$I$,
our extension is somewhat limited.

\begin{thm}\label{thm1}
Let $M$ be a rational polygonal right prism of the form \eqref{eq1.1}, where $P$ is a rational polygon and $I=[0,z_0]$ is an interval.
For almost every pair of initial direction and starting point, the half-infinite billiard orbit in $M$ is uniformly distributed.
\end{thm}

For illustration, we consider a special case where $P$ is a right triangle.
It is well known that the right triangle billiard with angle $\pi/4$ and the right triangle billiard with angle $\pi/6$
are the only right triangle billiards that are integrable, exhibiting stable and predictable behaviour.
Perhaps the simplest non-integrable billiard is the right triangle billiard with angle~$\pi/8$.
It is also well known that unfolding in the spirit of K\"{o}nig and Sz\"{u}cs~\cite{KS13} leads to a $16$-fold covering
of the triangle and shows that this billiard is equivalent to geodesic flow on the regular octagon surface $\PPP$
where parallel edges are identified in pairs.
On the other hand, unfolding also leads to a $2$-fold covering of the interval $I=[0,z_0]$.
Thus billiard in the rational polygonal right prism $M=P\times I$, where $P$ is the right triangle with angle $\pi/8$
and $I=[0,z_0]$, is equivalent to geodesic flow in the translation $3$-manifold $\MMM=\PPP\times\III$,
where $\PPP$ is the regular octagon translation surface and $\III=[0,2z_0]$, treated as a torus.
Figure~1 illustrates that $\MMM=\PPP\times\III$ gives a $32$-fold covering of the rational polygonal right prism $M=P\times I$.
It has $2$ octagonal faces which are identified with each other, and $8$ rectangular faces, with pairs of parallel ones
identified with each other, analogous to the edge identification of the regular octagon translation surface~$\PPP$.

\begin{displaymath}
\begin{array}{c}
\includegraphics[scale=0.8]{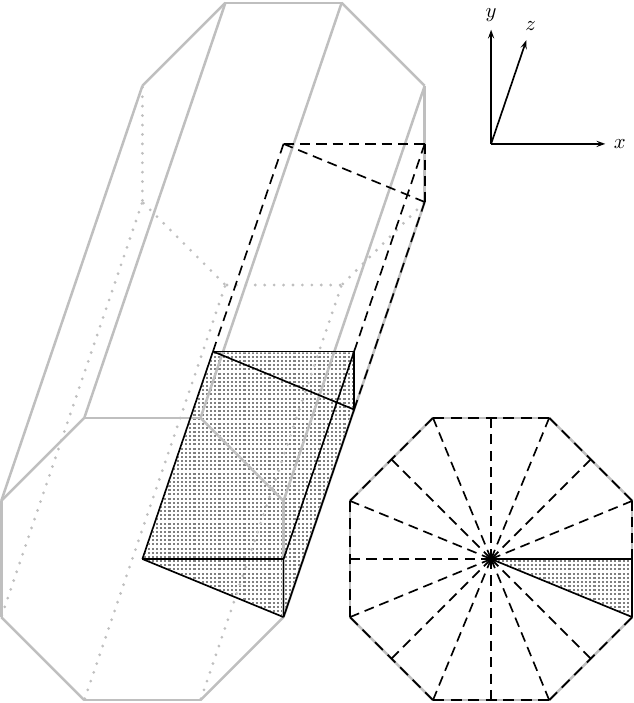}
\vspace{3pt}\\
\mbox{Figure 1: the translation $3$-manifold $\MMM=\PPP\times\III$}
\end{array}
\end{displaymath}

\begin{thm}\label{thm2}
Let $\MMM=\PPP\times\III$ be a rational octagonal right prism translation $3$-manifold,
where $\PPP$ is the regular octagon translation surface and $\III=[0,2z_0]$, treated as a torus.
For almost every pair of direction and starting point, the half-infinite geodesic in $\MMM$ is uniformly distributed.
\end{thm}

%
%

\section{Proof of Theorem~\ref{thm2}}\label{sec2}

Suppose that a half-infinite geodesic
\begin{displaymath}
\LLL_{S_0,\bfv}(t)=(s_1+v_1t,s_2+v_2t,s_3+v_3t),
\quad
t\ge0,
\end{displaymath}
in $\MMM$ has a non-pathological starting point $\LLL_{S_0,\bfv}(0)=S_0=(s_1,s_2,s_3)$,
and direction given by the unit vector
\begin{displaymath}
\bfv=(v_1,v_2,v_3)\in\Rr^3,
\quad\mbox{where}\quad
v_1^2+v_2^2+v_3^2=1,
\end{displaymath}
with arc-length parametrization.
The coordinates $(s_1+v_1t,s_2+v_2t)$ are modulo~$\PPP$ and the coordinate $s_3+v_3t$ is modulo~$\III$.

We may assume without loss of generality that $v_3>0$.
Then the geodesic $\LLL_{S_0,\bfv}(t)$, $t\ge0$, hits the octagon face of $\MMM$ for the very first time at time $t=t_0$,
where $s_3+v_3t_0=2z_0$, so that $t_0=(2z_0-s_3)/v_3$.
Indeed, the geodesic hits the octagon face of $\MMM$ for the $(k+1)$-th time at time $t=t_k$, where
\begin{equation}\label{eq2.1}
t_k=\frac{2kz_0+(2z_0-s_3)}{v_3}=k\theta+\lambda,
\quad
k=0,1,2,3,\ldots,
\end{equation}
with parameters
\textcolor{white}{xxxxxxxxxxxxxxxxxxxxxxxxxxxxxx}
\begin{displaymath}
\theta=\frac{2z_0}{v_3}\quad\mbox{and}\quad
\lambda=\frac{2z_0-s_3}{v_3}.
\end{displaymath}
This gives rise to an arithmetic progression
\begin{displaymath}
\lambda<\theta+\lambda<2\theta+\lambda<3\theta+\lambda<\ldots,
\end{displaymath}
with common gap $\theta$ between consecutive terms.

We need the following result on arithmetic progressions; see \cite[Theorem~2.2.2]{BDY20a}.

\begin{lem}[transference theorem for arithmetic progressions]\label{lem21}
Let $\SSS\subset\Rr$ be a measurable set.
For every $\ell\in\Zz$, there exists a constant $c_1(\ell)>0$, dependent only on~$\ell$, such that
for almost every pair $\theta,\lambda$ satisfying $2^\ell\le\theta<2^{\ell+1}$ and $0\le\lambda<\theta$,
the inequality
\begin{displaymath}
\left\vert\sum_{\substack{{k\ge0}\\{k\theta+\lambda\in\SSS\cap[0,n]}}}1-\frac{1}{\theta}\meas(\SSS\cap[0,n])\right\vert
\le c_1(\ell)n^{3/4}(\log n)^{1/2}
\end{displaymath}
holds for every sufficiently large positive integer~$n$.
\end{lem}

We have following immediate consequence.

\begin{lem}\label{lem22}
Suppose that the set $\SSS\subset\Rr$ is measurable.
Suppose further that $\SSS$ has asymptotic density $d=d(\SSS)\in[0,1]$,
so that there exists a monotonic sequence $\eps(n)=\eps(\SSS;n)\to0$ as $n\to\infty$ such that
\begin{displaymath}
\left\vert\frac{1}{n}\meas(\SSS\cap[0,n])-d(\SSS)\right\vert<\eps(n),
\quad
n=1,2,3,\ldots.
\end{displaymath}
For every $\ell\in\Zz$, there exists a constant $c_2(\ell)>0$, dependent only on~$\ell$, such that
for almost every pair $\theta,\lambda$ satisfying $2^\ell\le\theta<2^{\ell+1}$ and $0\le\lambda<\theta$,
the inequality
\begin{equation}\label{eq2.2}
\left\vert\frac{\theta}{n}\sum_{\substack{{k\ge0}\\{k\theta+\lambda\in\SSS\cap[0,n]}}}1-d(\SSS)\right\vert
\le\eps(n)+c_2(\ell)\frac{(\log n)^{1/2}}{n^{1/4}}
\end{equation}
holds for every sufficiently large positive integer~$n$.
\end{lem}

Meanwhile, every point in $\MMM$ is of the form $(x,y,z)$, where $(x,y)\in\PPP$ and $z\in\III$.
We consider the projection
\begin{equation}\label{eq2.3}
\phi:\MMM\to\PPP:(x,y,z)\mapsto(x,y).
\end{equation}
Then the image of the geodesic $\LLL_{S_0,\bfv}(t)$, $t\ge0$, under this projection is a half-infinite geodesic
\textcolor{white}{xxxxxxxxxxxxxxxxxxxxxxxxxxxxxx}
\begin{equation}\label{eq2.4}
\HHH_{S_0,\bfv}(t)=(s_1+v_1t,s_2+v_2t),
\quad
t\ge0,
\end{equation}
on the regular octagon translation surface~$\PPP$.
Clearly the key parameters $s_3$ and~$v_3$, particularly concerning the hitting times given in \eqref{eq2.1},
are lost under this projection \eqref{eq2.3}.
However, we know the arithmetic progression \eqref{eq2.1} of the time instances when the geodesic $\LLL_{S_0,\bfv}(t)$, $t\ge0$,
hits the octagon face of~$\MMM$.
This gives rise to an infinite sequence of points $\HHH_{S_0,\bfv}(t_k)$, $k=0,1,2,3,\ldots,$ on~$\PPP$.
For any $S_0$ and~$\bfv$, if we can show that this sequence of points is uniformly distributed on~$\PPP$,
then the half-infinite geodesic $\LLL_{S_0,\bfv}(t)$, $t\ge0$, is uniformly distributed in~$\MMM$.

Consider a typical half-infinite geodesic $\HHH_\bfw(\tau)$, $\tau\ge0$, on~$\PPP$,
with direction given by the unit vector $\bfw=(w_1,w_2)\in\Rr^2$.
Suppose that this geodesic is the image on $\PPP$ of $\LLL_{S_0,\bfv}(t)$, $t\ge0$, under the projection \eqref{eq2.3}.
Then
\begin{equation}\label{eq2.5}
\HHH_\bfw(\tau)=(s_1+w_1\tau,s_2+w_2\tau),
\quad
\tau\ge0.
\end{equation}
In view of the different parametrizations of \eqref{eq2.4} and \eqref{eq2.5}, we have
\begin{displaymath}
\HHH_{S_0,\bfv}(t)=\HHH_\bfw(\tau)
\quad\mbox{if and only if}\quad
\tau=(v_1^2+v_2^2)^{1/2}t.
\end{displaymath}
Corresponding to the arithmetic progression $t_k$, $k=0,1,2,3,\ldots,$ of hitting times given by \eqref{eq2.1} is the arithmetic progression
\begin{displaymath}
\tau_k=(v_1^2+v_2^2)^{1/2}t_k=\frac{(v_1^2+v_2^2)^{1/2}(2kz_0-s_3)}{v_3}=k\theta+\lambda,
\quad
k=0,1,2,3,\ldots,
\end{displaymath}
with parameters
\textcolor{white}{xxxxxxxxxxxxxxxxxxxxxxxxxxxxxx}
\begin{equation}\label{eq2.6}
\theta=\frac{2z_0(v_1^2+v_2^2)^{1/2}}{v_3}\quad\mbox{and}\quad
\lambda=\frac{(2z_0-s_3)(v_1^2+v_2^2)^{1/2}}{v_3}.
\end{equation}

By the geodesic analogue of Theorem~A, for almost every direction $\bfw$ and for every non-pathological starting point for this direction,
the geodesic $\HHH_\bfw(\tau)$, $\tau\ge0$, is uniformly distributed on~$\PPP$.
Let $\rp$ denote an arbitrary polygon on $\PPP$ where all the vertices have rational coordinates, and let the measurable set
\begin{equation}\label{eq2.7}
\SSS=\SSS(\HHH_\bfw;\rp)=\{\tau\ge0:\HHH_\bfw(\tau)\in\rp\}
\end{equation}
denote the set of time instances when this geodesic visits~$\rp$.
The uniformity of the geodesic then implies that $\SSS$ has asymptotic density
\begin{displaymath}
d=d(\SSS)=\frac{\area(\rp)}{\area(\PPP)}\in[0,1].
\end{displaymath}

The uniformly distributed geodesic $\HHH_\bfw(\tau)$, $\tau\ge0$, is clearly the image on $\PPP$ under the projection \eqref{eq2.3}
of infinitely many different geodesics $\LLL_{S_0,\bfv}(t)$, $t\ge0$, in the $3$-manifold~$\MMM$, as there are only two requirements,
namely $\HHH_\bfw(0)=(s_1,s_2)$ concerning the starting point, where $S_0=(s_1,s_2,s_3)$,
and $w_1v_2=w_2v_1$ concerning equality of the relevant directions.
Applying Lemma~\ref{lem22}, we see that for every $\ell\in\Zz$ and for almost every pair $\theta,\lambda$ of the form \eqref{eq2.6}
satisfying $2^\ell\le\theta<2^{\ell+1}$ and $0\le\lambda<\theta$, the inequality \eqref{eq2.2} with $\SSS$ given by \eqref{eq2.7}
holds for every sufficiently large positive integer~$n$.

This means that for almost every pair of starting point $S_0$ and unit direction vector~$\bfv$, the infinite sequence
\begin{equation}\label{eq2.8}
\LLL_{S_0;\bfv}(t_k),
\quad
k=0,1,2,3,\ldots,
\end{equation}
of points, where the sequence $t_k$, $k=0,1,2,3,\ldots,$ of time instances is given by \eqref{eq2.1}, is uniformly distributed on $\PPP$
relative to the single test set~$\rp$.

The set of all polygons $\rp$ on $\PPP$ where all the vertices have rational coordinates is countable.
On the other hand, a countable union of sets of measure zero has measure zero.
It follows that for almost every pair of starting point $S_0$ and unit direction vector~$\bfv$, the infinite sequence \eqref{eq2.8}
of points is uniformly distributed on $\PPP$ relative to every polygon $\rp$ on $\PPP$ where all the vertices have rational coordinates.
This guarantees uniformity in general, in the classical Weyl sense, and completes the proof of Theorem~\ref{thm2}.

\begin{remark}
Theorem~\ref{thm1} is a result on time-qualitative uniformity, and does not say anything about the speed of convergence
to uniform distribution, as a key ingredient of the proof is the geodesic analogue of Theorem~A which is also time-qualitative in nature.
There are instances, however, when we can establish time-quantitative results.
This happens, for example, when it is possible to establish time-quantitative uniform distribution results for some geodesics on the underlying
rational polygonal translation surface~$\PPP$.
We can establish such extensions of the geodesic analogue of Theorem~A in \cite{BDY20a,BDY20b} for the L-surface
and in \cite{BCY22b} for finite polysquare translation surfaces and for regular polygonal translation surfaces.
These in turn lead to extensions of various analogues of Theorem~\ref{thm2}, and hence also Theorem~\ref{thm1},
to time-quantitative results.
\end{remark}

%
%

\section{Proof of Lemma~\ref{lem21}}\label{sec3}

Throughout the proof, the set $\SSS\subset\Rr$ is measurable.

Let a non-negative integer $\ell$ be chosen and fixed.

Consider an infinite sequence $N_1,N_2,N_3,\ldots$ of positive integers satisfying
\begin{displaymath}
1<N_1<N_2<N_3<\ldots<N_h<\ldots,
\end{displaymath}
and another infinite sequence $M_1,M_2,M_3,\ldots$ of positive integers satisfying
\begin{displaymath}
1<M_h<N_h,
\quad
h=1,2,3,\ldots,
\end{displaymath}
both to be specified later in terms of the parameter $h$ and the chosen integer~$\ell$.
For every positive integer~$h$, let $S(h)\subset[0,1]$ denote the \textit{contraction} of $\SSS\cap[0,N_h]$
to the unit interval, so that
\begin{equation}\label{eq3.1}
x\in S(h)
\quad\mbox{if and only if}\quad
N_hx\in\SSS\cap[0,N_h].
\end{equation}
Since the characteristic function
\begin{displaymath}
\chi_{S(h)}(x)=\left\{\begin{array}{ll}
1,&\mbox{if $x\in S(h)$},\\
0,&\mbox{if $x\not\in S(h)$},
\end{array}\right.
\end{displaymath}
defined over $[0,1]$ and extended periodically over the whole real line with period~$1$, is measurable,
we can consider its Fourier series
\begin{equation}\label{eq3.2}
\chi_{S(h)}(x)=\sum_{j\in\Zz}a_j\ee^{2\pi\ii jx},
\end{equation}
with Fourier coefficients $a_j$, $j\in\Zz$.
In particular,
\begin{displaymath}
a_0=\lambda_1(S(h)).
\end{displaymath}

\begin{remark}
For a measurable set $S(h)$, the infinite Fourier series \eqref{eq3.2} may diverge at some points.
However, Lemma~\ref{lem21} is a measure theoretic statement which ignores sets of measure zero.
So it suffices to have pointwise convergence almost everywhere.
Fourier analysis provides at least two options to settle this issue.
We can use the very deep Carleson theorem.
Alternatively, we can use the much simpler Lebesgue theorem with Ces\`{a}ro summability.
\end{remark}

The Parseval formula gives
\begin{displaymath}
\sum_{j\in\Zz}\vert a_j\vert^2=\lambda_1(S(h))\le1,
\end{displaymath}
so that
\begin{equation}\label{eq3.3}
\sum_{j\in\Zz\setminus\{0\}}\vert a_j\vert^2
=\lambda_1(S(h))-\lambda_1^2(S(h))
=\lambda_1(S(h))(1-\lambda_1(S(h)))<1.
\end{equation}

Lemma~\ref{lem21} for the chosen value $\ell$ concerns the arithmetic progression $k\theta+\eta$, $k\ge0$,
where $2^\ell\le\theta<2^{\ell+1}$ and $0\le\eta<\theta$.
The contraction \eqref{eq3.1} leads to a new arithmetic progression $k\alpha+\beta$, $k\ge0$,
where $2^\ell/N_h\le\alpha<2^{\ell+1}/N_h$ and $0\le\beta<\alpha$.

For any $\alpha\in[2^\ell/N_h,2^{\ell+1}/N_h)$, let $K(\alpha)$ be the unique integer satisfying
\begin{equation}\label{eq3.4}
(K(\alpha)-1)\alpha<1\le K(\alpha)\alpha.
\end{equation}
Using the Fourier series \eqref{eq3.2}, we have
\begin{equation}\label{eq3.5}
\sum_{k=0}^{K(\alpha)-1}\chi_{S(h)}(k\alpha+\beta)-K(\alpha)\lambda_1(S(h))
=\sum_{j\in\Zz\setminus\{0\}}a_j\sum_{k=0}^{K(\alpha)-1}\ee^{2\pi\ii j(k\alpha+\beta)}
\end{equation}
for every $\alpha$ and $\beta$ satisfying $2^\ell/N_h\le\alpha<2^{\ell+1}/N_h$ and $0\le\beta<\alpha$.

To study \eqref{eq3.5}, we consider the integral
\begin{equation}\label{eq3.6}
J(\bfa;N_h;M_h)=\int_{2^\ell/N_h}^{2^{\ell+1}/N_h}\!\!\int_{-2^{\ell-1}/M_h}^{2^{\ell-1}/M_h}
\left\vert\sum_{j\in\Zz\setminus\{0\}}a_j\sum_{k=0}^{K(\alpha)-1}\ee^{2\pi\ii j(k\alpha+\gamma)}\right\vert^2
\dd\gamma\,\dd\alpha.
\end{equation}
To obtain a bound on this integral, we observe that for $\gamma\in[-2^{\ell-1}/M_h,2^{\ell-1}/M_h]$, the inequality
$2(1-\vert\gamma\vert M_h/2^\ell)\ge1$ holds.
It then follows that
\begin{equation}\label{eq3.7}
J(\bfa;N_h;M_h)\le2J^*(\bfa;N_h;M_h),
\end{equation}
where for integers $N$ and $M$ satisfying $1<M<N$,
\begin{align}\label{eq3.8}
&
J^*(\bfa;N;M)
\nonumber
\\
&\quad
=\int_{2^\ell/N}^{2^{\ell+1}/N}\!\!\int_{-2^\ell/M}^{2^\ell/M}
\left\vert\sum_{j\in\Zz\setminus\{0\}}a_j\sum_{k=0}^{K(\alpha)-1}\ee^{2\pi\ii j(k\alpha+\gamma)}\right\vert^2
\left(1-\frac{\vert\gamma\vert M}{2^\ell}\right)\,\dd\gamma\,\dd\alpha.
\end{align}
At the end of this section, we establish the following bound for this integral.

\begin{lem}\label{lem31}
For any sequence $\bfa$ satisfying \eqref{eq3.3}, the inequality
\begin{equation}\label{eq3.9}
J^*(\bfa;N;M)\le2^{\ell+11}.
\end{equation}
holds uniformly for integers $M$ and $N$ satisfying $1<M<N$,
where $K(\alpha)$ is the integer defined by \eqref{eq3.4} for every $\alpha\in[2^\ell/N,2^{\ell+1}/N)$.
\end{lem}

Combining \eqref{eq3.6}, \eqref{eq3.7} and \eqref{eq3.9}, we deduce that
\begin{equation}\label{eq3.10}
J(\bfa;N_h;M_h)\le2^{\ell+12}.
\end{equation}
The inequality \eqref{eq3.10} is a quadratic average result, from which we can derive information
concerning the majority of pairs
\begin{equation}\label{eq3.11}
(\alpha,\gamma)\in\left[\frac{2^\ell}{N_h},\frac{2^{\ell+1}}{N_h}\right)
\times\left[-\frac{2^{\ell-1}}{M_h},\frac{2^{\ell-1}}{M_h}\right].
\end{equation}
Let
\textcolor{white}{xxxxxxxxxxxxxxxxxxxxxxxxxxxxxx}
\begin{displaymath}
B(N_h;M_h)=\left\{(\alpha,\gamma)\in\left[\frac{2^\ell}{N_h},\frac{2^{\ell+1}}{N_h}\right)
\times\left[-\frac{2^{\ell-1}}{M_h},\frac{2^{\ell-1}}{M_h}\right]:
\mbox{\eqref{eq3.12} holds}\right\}
\end{displaymath}
denote the collection of pairs $(\alpha,\gamma)$ satisfying \eqref{eq3.11} such that
\begin{equation}\label{eq3.12}
\left\vert\sum_{j\in\Zz\setminus\{0\}}a_j\sum_{k=0}^{K(\alpha)-1}\ee^{2\pi\ii j(k\alpha+\gamma)}\right\vert
\ge(hN_hM_h)^{1/2}\log(1+h).
\end{equation}
Then it follows from \eqref{eq3.10} that
\begin{equation}\label{eq3.13}
\lambda_2(B(N_h;M_h))\le\frac{2^{\ell+12}}{hN_hM_h\log^2(1+h)},
\end{equation}
where $\lambda_2$ denotes $2$-dimensional Lebesgue measure.

Next, note that as we move from \eqref{eq3.5} to \eqref{eq3.6},
we replace the parameter $\beta$ over a short interval $[0,\alpha)$
by a parameter $\gamma\in[-2^{\ell-1}/M_h,2^{\ell-1}/M_h]$ over a longer interval.
For any pair
\textcolor{white}{xxxxxxxxxxxxxxxxxxxxxxxxxxxxxx}
\begin{equation}\label{eq3.14}
(\alpha,\beta)\in\left[\frac{2^\ell}{N_h},\frac{2^{\ell+1}}{N_h}\right)\times[0,\alpha),
\end{equation}
there are at least $N_h/2M_h$ values of $\gamma\in[-2^{\ell-1}/M_h,2^{\ell-1}/M_h]$ where $\{\gamma/\alpha\}\alpha=\beta$.
For each of these values of~$\gamma$, consider the two arithmetic progressions
\begin{equation}\label{eq3.15}
k\alpha+\gamma,
\quad
k=0,1,2,3,\ldots,K(\alpha)-1,
\end{equation}
and
\textcolor{white}{xxxxxxxxxxxxxxxxxxxxxxxxxxxxxx}
\begin{equation}\label{eq3.16}
k\alpha+\beta,
\quad
k=0,1,2,3,\ldots,K(\alpha)-1.
\end{equation}
Since
\textcolor{white}{xxxxxxxxxxxxxxxxxxxxxxxxxxxxxx}
\begin{displaymath}
\gamma=\left[\frac{\gamma}{\alpha}\right]\alpha+\beta,
\end{displaymath}
the arithmetic progression \eqref{eq3.15} is obtained
by simply advancing the arithmetic progression \eqref{eq3.16} by $[\gamma/\alpha]$ terms.
More precisely, the arithmetic progression \eqref{eq3.15} is given by
\begin{equation}\label{eq3.17}
k\alpha+\beta,
\quad
k=\left[\frac{\gamma}{\alpha}\right],\left[\frac{\gamma}{\alpha}\right]+1,\left[\frac{\gamma}{\alpha}\right]+2,\left[\frac{\gamma}{\alpha}\right]+3,
\ldots,\left[\frac{\gamma}{\alpha}\right]+K(\alpha)-1.
\end{equation}

\begin{lem}\label{lem32}
If a pair $(\alpha,\beta)$ such that \eqref{eq3.14} holds satisfies the inequality
\begin{equation}\label{eq3.18}
\left\vert\sum_{j\in\Zz\setminus\{0\}}a_j\sum_{k=0}^{K(\alpha)-1}\ee^{2\pi\ii j(k\alpha+\beta)}\right\vert
\ge\frac{2N_h}{M_h}+(hN_hM_h)^{1/2}\log(1+h),
\end{equation}
then each pair $(\alpha,\gamma)$ such that \eqref{eq3.11} and $\{\gamma/\alpha\}\alpha=\beta$ hold
satisfies the inequality \eqref{eq3.12}.
\end{lem}

\begin{proof}
It clearly suffices to prove that
\begin{equation}\label{eq3.19}
\left\vert
\sum_{j\in\Zz\setminus\{0\}}a_j\sum_{k=0}^{K(\alpha)-1}\ee^{2\pi\ii j(k\alpha+\beta)}
-\sum_{j\in\Zz\setminus\{0\}}a_j\sum_{k=0}^{K(\alpha)-1}\ee^{2\pi\ii j(k\alpha+\gamma)}
\right\vert
\le\frac{2N_h}{M_h}.
\end{equation}
Since
\textcolor{white}{xxxxxxxxxxxxxxxxxxxxxxxxxxxxxx}
\begin{displaymath}
\left\vert\left[\frac{\gamma}{\alpha}\right]\right\vert\le\frac{N_h}{2M_h},
\end{displaymath}
it follows from \eqref{eq3.16} and \eqref{eq3.17} that those terms
that belong to one of the arithmetic progressions \eqref{eq3.15} or \eqref{eq3.16}
but not both then form two arithmetic progressions of the form
\textcolor{white}{xxxxxxxxxxxxxxxxxxxxxxxxxxxxxx}
\begin{displaymath}
k\alpha+\rho,
\quad
k=0,1,2,3,\ldots,K-1,
\end{displaymath}
where $K\le N_h/2M_h$.
Hence
\begin{displaymath}
\sum_{k=0}^{K(\alpha)-1}\ee^{2\pi\ii j(k\alpha+\beta)}
-\sum_{k=0}^{K(\alpha)-1}\ee^{2\pi\ii j(k\alpha+\gamma)}
\end{displaymath}
is the sum of two sums of the form
\begin{displaymath}
\sum_{k=0}^{K-1}\ee^{2\pi\ii j(k\alpha+\rho)},
\end{displaymath}
where $K\le N_h/2M_h$.
Now for each of the two sums, we have
\begin{align}
\left\vert\sum_{j\in\Zz\setminus\{0\}}a_j\sum_{k=0}^{K-1}\ee^{2\pi\ii j(k\alpha+\rho)}\right\vert
&
=\left\vert\sum_{k=0}^{K-1}\sum_{j\in\Zz\setminus\{0\}}a_j\ee^{2\pi\ii j(k\alpha+\rho)}\right\vert
\nonumber
\\
&
=\left\vert\sum_{k=0}^{K-1}\left(\sum_{j\in\Zz}a_j\ee^{2\pi\ii j(k\alpha+\rho)}-a_0\right)\right\vert
\nonumber
\\
&
\le\sum_{k=0}^{K-1}\left\vert\chi_{S(h)}(k\alpha+\rho)-\lambda_1(S_h)\right\vert
\nonumber
\\
&
\le2K.
\nonumber
\end{align}
This clearly leads to \eqref{eq3.19} and completes the proof.
\end{proof}

Let
\textcolor{white}{xxxxxxxxxxxxxxxxxxxxxxxxxxxxxx}
\begin{displaymath}
A(N_h;M_h)=\left\{(\alpha,\beta)\in\left[\frac{2^\ell}{N_h},\frac{2^{\ell+1}}{N_h}\right]\times[0,\alpha):
\mbox{\eqref{eq3.18} holds}\right\}.
\end{displaymath}
Then the above argument leads to the inequality
\begin{equation}\label{eq3.20}
\lambda_2(B(N_h;M_h))\ge\frac{N_h}{2M_h}\lambda_2(A(N_h;M_h)).
\end{equation}
Combining \eqref{eq3.13} and \eqref{eq3.20}, we obtain the upper bound
\begin{displaymath}
\lambda_2(A(N_h;M_h))\le\frac{2^{\ell+13}}{hN_h^2\log^2(1+h)}.
\end{displaymath}
Combining this with \eqref{eq3.5}, it is not difficult to see that apart from a set of measure $\lambda_2(A(N_h;M_h))$,
every pair $(\alpha,\beta)$ such that \eqref{eq3.14} holds satisfies the inequality
\begin{displaymath}
\left\vert\sum_{k=0}^{K(\alpha)-1}\chi_{S(h)}(k\alpha+\beta)-K(\alpha)\lambda_1(S(h))\right\vert
<\frac{2N_h}{M_h}+(hN_hM_h)^{1/2}\log(1+h).
\end{displaymath}

Next, note that the two expressions
\begin{displaymath}
\sum_{\substack{{k\ge0}\\{k\alpha+\beta\in S(h)}}}1-\frac{1}{\alpha}\lambda_1(S(h))
\quad\mbox{and}\quad
\sum_{k=0}^{K(\alpha)-1}\chi_{S(h)}(k\alpha+\beta)-K(\alpha)\lambda_1(S(h))
\end{displaymath}
differ by at most~$2$, due to the possibility that $(K(\alpha)-1)\alpha+\beta>1$ and the difference
$\vert K(\alpha)-1/\alpha\vert<1$, in view of \eqref{eq3.4}.
It follows that on reversing the contraction, we see that apart from a set of measure at most
\begin{displaymath}
\frac{2^{\ell+13}}{h\log^2(1+h)},
\end{displaymath}
every pair $(\theta,\eta)$ such that $2^\ell\le\theta\le2^{\ell+1}$ and $\eta\in[0,\theta)$ satisfies the inequality
\begin{equation}\label{eq3.21}
\left\vert\sum_{\substack{{k\ge0}\\{k\theta+\eta\in\SSS\cap[0,N_h]}}}1-\frac{1}{\theta}\lambda_1(\SSS\cap[0,N_h])\right\vert
<\frac{2N_h}{M_h}+(hN_hM_h)^{1/2}\log(1+h)+2.
\end{equation}

\begin{lem}[Borel--Cantelli lemma]\label{lem33}
Let\index{Borel--Cantelli lemma}
$(X,\Sigma,\mu)$ be a measure space,\index{measure!space}
and suppose that $E_h$, $h=1,2,3,\ldots,$ is a sequence of $\Sigma$-measurable sets.
If
\begin{displaymath}
\sum_{h=1}^\infty\mu(E_h)<\infty,
\end{displaymath}
then
\textcolor{white}{xxxxxxxxxxxxxxxxxxxxxxxxxxxxxx}
\begin{displaymath}
\mu\left(\bigcap_{h=1}^\infty\bigcup_{i=h}^\infty E_i\right)=0.
\end{displaymath}
\end{lem}

Since
\textcolor{white}{xxxxxxxxxxxxxxxxxxxxxxxxxxxxxx}
\begin{displaymath}
\sum_{h=1}^\infty\frac{1}{h\log^2(1+h)}<\infty,
\end{displaymath}
we conclude that for almost every pair $(\theta,\eta)$ such that $2^\ell\le\theta\le2^{\ell+1}$ and $\eta\in[0,\theta)$,
the inequality \eqref{eq3.21} holds for all sufficiently large positive integers~$h$.

Finally, we specify the integers $N_h$ and $M_h$ in terms of the parameter $h\ge1$ and the chosen integer~$\ell$.
Choosing them to satisfy
\begin{equation}\label{eq3.22}
N_h\le2^\ell(h^4\log^2(1+h))<N_h+1
\quad\mbox{and}\quad
M_h=2^\ell h
\end{equation}
ensures that the two dominant terms on the right hand side of \eqref{eq3.21} have the same order of magnitude in terms of~$h$.
For an arbitrary sufficiently integer~$n$, we choose $h$ to satisfy $N_h\le n<N_{h+1}$.
Then it follows from \eqref{eq3.21} and \eqref{eq3.22} that
\begin{align}
&
\left\vert\sum_{\substack{{k\ge0}\\{k\theta+\eta\in\SSS\cap[0,n]}}}1-\frac{1}{\theta}\lambda_1(\SSS\cap[0,n])\right\vert
\nonumber
\\
&\qquad
<\frac{N_{h+1}-N_h}{\theta}+\frac{2N_h}{M_h}+(hN_hM_h)^{1/2}\log(1+h)+2
\nonumber
\\
&\qquad
\le\frac{2^\ell((h+1)^4\log^2(2+h)-h^4\log^2(1+h))}{\theta}+(2+2^\ell)h^3\log^2(1+h)+3
\nonumber
\\
&\qquad
\le2^{\ell+3}h^3\log^2(1+h)+O_\ell(h^3\log(1+h))
\le c_1(\ell)n^{3/4}(\log n)^{1/2},
\nonumber
\end{align}
provided that~$n$, and hence also~$h$, is sufficiently large.

This completes the proof of Lemma~\ref{lem21}.

\begin{proof}[Proof of Lemma~\ref{lem31}]
For any fixed $\delta\in(0,1/2)$, we define the roof function $R_\delta:\Rr\to\Rr$ by writing
\textcolor{white}{xxxxxxxxxxxxxxxxxxxxxxxxxxxxxx}
\begin{displaymath}
R_\delta(x)=\left\{\begin{array}{ll}
0,&\mbox{if $\vert x\vert>\delta$},\\
1-(\vert x\vert/\delta),&\mbox{if $0\le\vert x\vert\le\delta$}.
\end{array}\right.
\end{displaymath}
For every integer $j\in\Zz$, we consider the integral
\begin{equation}\label{eq3.23}
I(\delta;j)=\int_{-1/2}^{1/2}R_\delta(x)\ee^{2\pi\ii jx}\,\dd x.
\end{equation}
Then
\textcolor{white}{xxxxxxxxxxxxxxxxxxxxxxxxxxxxxx}
\begin{equation}\label{eq3.24}
I(\delta;0)=\delta
\quad\mbox{and}\quad
I(\delta;j)=\delta\left(\frac{\sin\pi j\delta}{\pi j\delta}\right)^2,
\quad
j\in\Zz\setminus\{0\}.
\end{equation}

For any integers $j_1,j_2\in\Zz\setminus\{0\}$ and positive integer~$N$, let
\begin{align}\label{eq3.25}
&
\frakB(j_1;j_2;N)
=\int_{2^\ell/N}^{2^{\ell+1}/N}
\left(a_{j_1}\sum_{k_1=0}^{K(\alpha)-1}\ee^{2\pi\ii j_1k_1\alpha}\right)
\left(\overline{a_{j_2}}\sum_{k_2=0}^{K(\alpha)-1}\ee^{-2\pi\ii j_2k_2\alpha}\right)\,\dd\alpha
\nonumber
\\
&\qquad
=\int_{2^\ell/N}^{2^{\ell+1}/N}
\left(a_{j_1}\frac{\ee^{2\pi\ii j_1K(\alpha)\alpha}-1}{\ee^{2\pi\ii j_1\alpha}-1}\right)
\left(\overline{a_{j_2}}\,\frac{\ee^{-2\pi\ii j_2K(\alpha)\alpha}-1}{\ee^{-2\pi\ii j_2\alpha}-1}\right)\dd\alpha,
\end{align}
so that
\begin{align}\label{eq3.26}
&
\vert \frakB(j_1;j_2;N)\vert
\nonumber
\\
&\quad
\le\frac{1}{2}\int_{2^\ell/N}^{2^{\ell+1}/N}
\left(\left\vert a_{j_1}\frac{\ee^{2\pi\ii j_1K(\alpha)\alpha}-1}{\ee^{2\pi\ii j_1\alpha}-1}\right\vert^2
+\left\vert a_{j_2}\frac{\ee^{2\pi\ii j_2K(\alpha)\alpha}-1}{\ee^{2\pi\ii j_2\alpha}-1}\right\vert^2\right)\dd\alpha.
\end{align}
Then it follows from \eqref{eq3.8}, \eqref{eq3.23} with $\delta=2^\ell/M$ and $j=j_1-j_2$
and from \eqref{eq3.25} that
\textcolor{white}{xxxxxxxxxxxxxxxxxxxxxxxxxxxxxx}
\begin{displaymath}
J^*(\bfa;N;M)
=\sum_{j_1\in\Zz\setminus\{0\}}\sum_{j_2\in\Zz\setminus\{0\}}
I\left(\frac{2^\ell}{M};j_1-j_2\right)\frakB(j_1;j_2;N).
\end{displaymath}
We write
\textcolor{white}{xxxxxxxxxxxxxxxxxxxxxxxxxxxxxx}
\begin{equation}\label{eq3.27}
J^*(\bfa;N;M)=J^*_1(\bfa;N;M)+J^*_2(\bfa;N;M),
\end{equation}
where $J^*_1(\bfa;N;M)$ contains all the diagonal terms in $J^*(\bfa;N;M)$ with $j_1=j_2$,
while $J^*_2(\bfa;N;M)$ contains all the off-diagonal terms in $J^*(\bfa;N;M)$ with $j_1\ne j_2$.
Noting that $I(2^\ell/M;0)=2^\ell/M$, we see that
\begin{equation}\label{eq3.28}
J^*_1(\bfa;N;M)=\frac{2^\ell}{M}\sum_{j\in\Zz\setminus\{0\}}\vert a_j\vert^2E(j;N),
\end{equation}
where
\textcolor{white}{xxxxxxxxxxxxxxxxxxxxxxxxxxxxxx}
\begin{equation}\label{eq3.29}
E(j;N)=\int_{2^\ell/N}^{2^{\ell+1}/N}\left\vert\frac{\ee^{2\pi\ii jK(\alpha)\alpha}-1}{\ee^{2\pi\ii j\alpha}-1}\right\vert^2\dd\alpha.
\end{equation}
Meanwhile, noting \eqref{eq3.24}, we see that
\begin{equation}\label{eq3.30}
J^*_2(\bfa;N;M)=\frac{2^\ell}{M}\mathop{\sum_{j_1\in\Zz\setminus\{0\}}\sum_{j_2\in\Zz\setminus\{0\}}}_{j_1\ne j_2}
\left(\frac{\sin\pi(j_1-j_2)2^\ell M^{-1}}{\pi(j_1-j_2)2^\ell M^{-1}}\right)^2\frakB(j_1;j_2;N).
\end{equation}
Combining \eqref{eq3.26} and \eqref{eq3.30}, we deduce that
\begin{equation}\label{eq3.31}
\vert J^*_2(\bfa;N;M)\vert
\le\frac{2^\ell}{M}\sum_{j\in\Zz\setminus\{0\}}\sum_{\zeta\in\Zz\setminus\{0\}}
\left(\frac{\sin\pi\zeta2^\ell M^{-1}}{\pi\zeta2^\ell M^{-1}}\right)^2
\vert a_j\vert^2E(j;N).
\end{equation}
It then follows from \eqref{eq3.27}, \eqref{eq3.28} and \eqref{eq3.31} that
\begin{equation}\label{eq3.32}
\!\!\!
\vert J^*(\bfa;N;M)\vert
\le\frac{2^\ell}{M}\left(1+\sum_{\zeta\in\Zz\setminus\{0\}}
\left(\frac{\sin\pi\zeta2^\ell M^{-1}}{\pi\zeta2^\ell M^{-1}}\right)^2\right)
\sum_{j\in\Zz\setminus\{0\}}\vert a_j\vert^2E(j;N).
\end{equation}
Next, note that
\begin{align}\label{eq3.33}
&
\frac{2^\ell}{M}\left(1+\sum_{\zeta\in\Zz\setminus\{0\}}
\left(\frac{\sin\pi\zeta2^\ell M^{-1}}{\pi\zeta2^\ell M^{-1}}\right)^2\right)
\le\frac{2^\ell}{M}\left(\sum_{\vert\zeta\vert\le M/2^\ell}1
+\sum_{\vert\zeta\vert>M/2^\ell}\left(\frac{M}{\pi\zeta2^\ell}\right)^2\right)
\nonumber
\\
&\qquad
\le\frac{2^\ell}{M}\left(\frac{2M}{2^\ell}+1
+\frac{2M^2}{\pi^24^\ell}\int_{M/2^\ell}^\infty\frac{\dd x}{x^2}\right)
\le6.
\end{align}

To complete the proof of Lemma~\ref{lem31}, in view of \eqref{eq3.3}, \eqref{eq3.32}
and \eqref{eq3.33}, it suffices to show that for every $j\in\Zz\setminus\{0\}$ and integer $N>1$, we have
\begin{equation}\label{eq3.34}
E(j;N)\le2^{\ell+8}.
\end{equation}

Consider first small values of~$j$, where $1\le\vert j\vert\le N/2^{\ell+2}$.
Suppose first that $j>0$.
Since $\alpha\in[2^\ell/N,2^{\ell+1}/N)$, it follows that $0\le j\alpha\le1/2$,
so that $\ee^{2\pi\ii j\alpha}$ is a point on the upper half circle of unit radius.
On the other hand, \eqref{eq3.4} implies the inequality $0\le jK(\alpha)\alpha-j<j\alpha$.
Hence $\ee^{2\pi\ii jK(\alpha)\alpha}=\ee^{2\pi\ii(jK(\alpha)\alpha-j)}$ is a point on the circular arc of
the upper half circle of unit radius joining the points $1$ and $\ee^{2\pi\ii j\alpha}$.
As shown in Figure~4.1, we clearly have $\vert\ee^{2\pi\ii jK(\alpha)\alpha}-1\vert<\vert\ee^{2\pi\ii j\alpha}-1\vert$.

\begin{displaymath}
\begin{array}{c}
\includegraphics[scale=0.8]{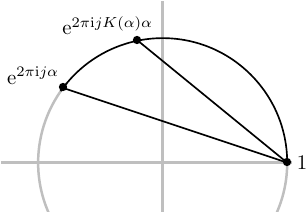}
\vspace{3pt}\\
\mbox{Figure 2: justifying the inequality \eqref{eq3.35}}
\end{array}
\end{displaymath}

Meanwhile, replacing $j$ by $-j$ preserves this inequality.
It follows that for every integer $j$ satisfying $1\le\vert j\vert\le N/2^{\ell+2}$, we have
\begin{equation}\label{eq3.35}
\left\vert\frac{\ee^{2\pi\ii jK(\alpha)\alpha}-1}{\ee^{2\pi\ii j\alpha}-1}\right\vert<1
\quad\mbox{and so}\quad
E(j;N)\le\frac{2^\ell}{N}.
\end{equation}

Next, for any fixed integer $j$ satisfying $\vert j\vert>N/2^{\ell+2}$,
we use the inequalities
\begin{equation}\label{eq3.36}
\left\vert\frac{\ee^{2\pi\ii jK(\alpha)\alpha}-1}{\ee^{2\pi\ii j\alpha}-1}\right\vert
\le\min\left\{K(\alpha),\frac{1}{\Vert j\alpha\Vert}\right\}
\le\min\left\{N,\frac{1}{\Vert j\alpha\Vert}\right\},
\end{equation}
where $\Vert x\Vert$ denotes the distance of a real number $x$ from the nearest integer,
and the inequality $K(\alpha)\le N$ follows from $\alpha\in[2^\ell/N,2^{\ell+1}/N)$ and \eqref{eq3.4}.

Let $r$ be the integer closest to~$j\alpha$.
Then
\begin{displaymath}
\Vert j\alpha\Vert<\frac{1}{N}
\quad\mbox{if and only if}\quad
\left\vert\alpha-\frac{r}{j}\right\vert<\frac{1}{\vert j\vert N},
\end{displaymath}
so that
\textcolor{white}{xxxxxxxxxxxxxxxxxxxxxxxxxxxxxx}
\begin{equation}\label{eq3.37}
\Vert j\alpha\Vert<\frac{1}{N}
\quad\mbox{if and only if}\quad
\alpha\in\left(\frac{r}{j}-\frac{1}{\vert j\vert N},\frac{r}{j}+\frac{1}{\vert j\vert N}\right).
\end{equation}
For this fixed integer $j$ satisfying $\vert j\vert>N/2^{\ell+2}$, there are at most
\begin{equation}\label{eq3.38}
\max\left\{5,\frac{2^{\ell+2}\vert j\vert}{N}\right\}\le\frac{2^{\ell+5}\vert j\vert}{N}
\end{equation}
integers $r$ such that
\begin{equation}\label{eq3.39}
\left(\frac{r}{j}-\frac{1}{\vert j\vert N},\frac{r}{j}+\frac{1}{\vert j\vert N}\right)\cap\left[\frac{2^\ell}{N},\frac{2^{\ell+1}}{N}\right)
\ne\emptyset.
\end{equation}
On the other hand, suppose that $n$ is a positive fixed integer satisfying $2^n\le N$.
Then analogous to \eqref{eq3.37}, we have
\begin{equation}\label{eq3.40}
\frac{2^{n-1}}{N}\le\Vert j\alpha\Vert\le\frac{2^n}{N}
\quad\mbox{if and only if}\quad
\alpha\in\frakI(j;N;n),
\end{equation}
where
\textcolor{white}{xxxxxxxxxxxxxxxxxxxxxxxxxxxxxx}
\begin{equation}\label{eq3.41}
\frakI(j;N;n)
=\left[\frac{r}{j}-\frac{2^n}{\vert j\vert N},\frac{r}{j}-\frac{2^{n-1}}{\vert j\vert N}\right]
\cup\left[\frac{r}{j}+\frac{2^{n-1}}{\vert j\vert N},\frac{r}{j}+\frac{2^n}{\vert j\vert N}\right],
\end{equation}
except when $2^n/N>1/2$, in which case we have the modification
\begin{equation}\label{eq3.42}
\frakI(j;N;n)
=\left[\frac{r}{j}-\frac{1}{2\vert j\vert},\frac{r}{j}-\frac{2^{n-1}}{\vert j\vert N}\right]
\cup\left[\frac{r}{j}+\frac{2^{n-1}}{\vert j\vert N},\frac{r}{j}+\frac{1}{2\vert j\vert}\right].
\end{equation}
Similarly, for this fixed integer $j$ satisfying $\vert j\vert>N/2^{\ell+2}$ and fixed positive integer $n$ satisfying $2^n\le N$,
there are at most \eqref{eq3.38} integers $r$ such that
\begin{equation}\label{eq3.43}
\frakI(j;N;n)\cap\left[\frac{2^\ell}{N},\frac{2^{\ell+1}}{N}\right)\ne\emptyset.
\end{equation}

Combining \eqref{eq3.29} and \eqref{eq3.36}--\eqref{eq3.43},
we see that for any fixed integer $j$ satisfying $\vert j\vert>N/2^{\ell+2}$,
we have
\begin{align}
&
E(j;N)
\le\int_{2^\ell/N}^{2^{\ell+1}/N}\left(\min\left\{N,\frac{1}{\Vert j\alpha\Vert}\right\}\right)^2\dd\alpha
\nonumber
\\
&\qquad
\le N^2\,\frac{2^{\ell+5}\vert j\vert}{N}\,\frac{2}{\vert j\vert N}
+\sum_{\substack{{n=1}\\{2^n\le N}}}^\infty\left(\frac{N}{2^{n-1}}\right)^2\frac{2^{\ell+5}\vert j\vert}{N}\,\frac{2^n}{\vert j\vert N}
\le2^{\ell+5}\left(2+4\sum_{n=1}^\infty\frac{1}{2^n}\right),
\nonumber
\end{align}
confirming the assertion \eqref{eq3.34} and completing the proof.
\end{proof}

%
%


\begin{thebibliography}{9}

\bibitem{BDY20a}
J. Beck, M. Donders, Y. Yang.
Quantitative behavior of non-integrable systems I.
\textit{Acta Math. Hungar.} \textbf{161} (2020), 66--184.

\bibitem{BDY20b}
J. Beck, M. Donders, Y. Yang.
Quantitative behavior of non-integrable systems II.
\textit{Acta Math. Hungar.} \textbf{162} (2020), 220--324.

\bibitem{BCY22b}
J. Beck, W.W.L. Chen, Y. Yang.
Quantitative behavior of non-integrable systems IV.
\textit{Acta Math. Hungar.} \textbf{167} (2022), 1--160.

\bibitem{KMS86}
S. Kerckhoff, H. Masur, J. Smillie.
Ergodicity of billiard flows and quadratic differentials.
\textit{Ann. of Math.} \textbf{124} (1986), 293--311.

\bibitem{KS13}
D. K\"{o}nig, A. Sz\"{u}cs.
Mouvement d'un point abondonne a l'interieur d'un cube.
\textit{Rend. Circ. Mat. Palermo} \textbf{36} (1913), 79--90.

\end{thebibliography}
\end{document}